\documentclass{amsart}
\usepackage{amsmath}
\usepackage[unicode,colorlinks=true,linkcolor=blue,urlcolor=magenta, citecolor=blue]{hyperref}

\usepackage[utf8x]{inputenc}
\usepackage{amssymb, upgreek}

\usepackage{color}
\definecolor{keywordcolor}{rgb}{0.7, 0.1, 0.1}   % red
\definecolor{commentcolor}{rgb}{0.4, 0.4, 0.4}   % grey
\definecolor{symbolcolor}{rgb}{0.0, 0.1, 0.6}    % blue
\definecolor{sortcolor}{rgb}{0.1, 0.5, 0.1}      % green
\definecolor{errorcolor}{rgb}{1, 0, 0}           % bright red
\definecolor{stringcolor}{rgb}{0.5, 0.3, 0.2}    % brown

\usepackage{listings}

\begin{document}

\title{What is the point of computers? A question for pure mathematicians}
\author{Kevin Buzzard}
\email{k.buzzard@imperial.ac.uk}
\address{Department of Mathematics, Imperial College London}
\begin{abstract}
  We discuss the idea that computers might soon help mathematicians to prove theorems in areas where they have not previously been useful. Furthermore we argue that these same computer tools will also help us in the communication and teaching of mathematics.
\end{abstract}

\maketitle

\section*{Introduction}

Computers in 2021 are phenomenal. They can do billions of calculations in a second. They are extremely good at obeying precise instructions accurately. Mathematics is a game with precise rules. One can thus ask in what ways computers can be used to help us\footnote{Throughout this article, by ``us'' and ``we'' I am referring to the community of people who, like myself, identify as pure mathematicians.} mathematicians to do our job.

%-- this is the science of mathematics. However a lot of human creativity is involved too, in deciding where we want to go and then figuring out how to get there -- this is the art of mathematics. too the art is in discovering how to apply them in the right order. 

Of course, computers have been used to help some mathematicians to do their job ever since computers have existed. Birch and Swinnerton-Dyer used an early computer (which was the size of a large room and which had 20 kilobytes of memory) to compute many examples of solutions to cubic equations in two variables modulo prime numbers~\cite{bsd}. Graphing the output data in the right way led to new insights in the theory of elliptic curves which ultimately became the Birch and Swinnerton-Dyer conjecture, one of the Clay Millennium problems.
% The conjecture states that two numerical invariants (the algebraic rank and the analytic rank) attached to an elliptic curve over the rational numbers are equal, and
At the time of writing, this conjecture is still open, although regular breakthroughs (most recently in non-commutative Iwasawa theory) provide us with incremental progress.

This article is not about using computers in that way. This article is an attempt to explain to \emph{all} researchers in mathematics that, thanks to breakthroughs in computer science, computers can now be used to help us not just with computations, but with \emph{reasoning}. In other words, it is about the possibility that computers might soon be helping us to \emph{prove theorems}, whether they be about ``computable'' objects such as elliptic curves, or about more intractable objects such as Banach spaces, schemes, abelian categories or perfectoid spaces, things which cannot be listed or classified or in general stored in a traditional computer algebra package in any meaningful way. In particular, it is about the possibility that computer proof assistants can help the mathematician who up until this point has had no need for computation in their research and might hence incorrectly deduce that computers have nothing at all to offer them. I should also stress that the applications are not limited to people interested in foundational subjects such as set theory or type theory; I am thinking about applications in geometry, topology, combinatorics, number theory, algebra, analysis,\ldots

I end this introduction with a summary of what to expect, and what not to expect, from this fast-growing area within the next decade. The first thing to stress is that computers will not be putting us out of a job. Computer proof assistants can now understand the \emph{statement} of the Riemann hypothesis, but I will eat my hat if a computer, all by itself, comes up with a \emph{proof} of the Riemann hypothesis (or indeed a proof of any open problem of interest to mainstream pure mathematicians) within the next 10 years\footnote{Conjectures which stretch beyond a 10 year period are I think very unwise; like mathematics, sometimes computer science moves very quickly.}.

What I do believe is going to happen within the next 10 years: tools will be created which will \emph{help} mathematicians to prove theorems. Digitised and semantically searchable databases of mathematics are appearing. Computers are going to start doing diagram chases for us, filling in the proofs of lemmas, pointing out counterexamples to our ideas, and suggesting results which might be helpful to us. The technology to make such tools is already coming; it is viable. Furthermore the databases of theorem statements and proofs which are appearing will not only have applications in research; we will be able to use them for teaching and for communicating mathematics in new ways. Undergraduates will be able to get instant feedback on their work. PhD students will be able to search for theorems and counterexamples in databases. Researchers will be able to write next-generation error-free papers where details can be folded and unfolded by the user. Patrick Massot has written a thoughtful piece~\cite{massot-note} explaining these and other ideas in more detail. Computers are going to be able to understand \emph{your area} of mathematics, and even keep up with it as it develops. But there is a catch. Who is going to make the database of important results in non-commutative Iwasawa theory, or whatever area you're interested in, which will power these tools? It's not going to be the computer scientists, because most of them know nothing about non-commutative Iwasawa theory. \emph{It has to be us.}

If you want to see progress within this domain in your own area of mathematics, I would \emph{urge} you to take some time working through some tutorials and learning one of these computer proof assistant languages. It is not difficult to do so -- I teach a popular course to final year mathematics students where we learn how to do undergraduate level mathematics (topology, analysis, group theory and so on) using the Lean theorem prover.\footnote{If you have Lean installed then you can take the course yourself; the materials are here~\cite{formalising-mathematics}.} Engaging with harder mathematics is not at all difficult \emph{once you know the language}. If you want to learn Lean's language, a good place to start is the Lean prover community's website~\cite{leanprover-community}. Coq and Isabelle/HOL are two other well-established theorem provers with big mathematics libraries, and there are plenty of others. If you can get to the point where you are able to explain the \emph{statements} of your own theorems to a computer proof assistant, then these statements can be added to databases, and furthermore you learnt a new skill. If however you can get to the point where you can explain the \emph{proofs}, then the AI people will be extremely interested, as will the people building huge formalised mathematical libraries which represent a 21st century Bourbaki. Furthermore, you will be having fun: formalisation of proofs is mathematics re-interpreted as an interesting computer puzzle game. If you don't have the time, then find a student who does. Instead of the traditional ``do a project consisting of reading a paper and then writing a paper showing that you understood the paper'', why not get a student to write some code which proves that they understood the paper? They can learn the language of the prover themselves, and then teach it to you as you teach them the mathematics.

The files which computer proof assistants can read and write represent a way of digitising mathematical ideas. Digitising something \emph{completely} changes (in fact it vastly augments) the ways in which it can be used. Consider for example the digitisation of music, with the CD and the mp3 file. This has revolutionised how music is consumed and delivered. My collection of music consists of hundreds of vinyl records, tapes and CDs in my office and loft. My children's collection is in the cloud, has essentially zero mass and volume, and is accessible anywhere. Not only that, but cloud based music platforms have also fundamentally changed the way the modern musician communicates with their fans, bypassing the traditional process completely. The music industry was turned upside-down by digitisation.

Mathematics has been done in the same ``pencil and paper'' way for millennia, but now there is a true opportunity to rethink and enhance this approach. I do not dare to dream what the ultimate consequences of digitising mathematics will be, but I firmly believe that it will make mathematics more accessible -- and easier for us to do, to communicate, and to play with. The ball is in our court.

\section*{Overview of the paper}

This paper describes a ``new'' way in which computers can be used by mathematicians. As mathematicians our typical experience with computers is that we can use traditional programming languages like python or traditional computer algebra packages like {\tt sage} to do things like compute the sum of the first 100 prime numbers. We know equally well that these traditional tools, even though they can compute as many prime numbers as you like (within reason), are not capable of \emph{proving} that there are infinitely many primes; the infinite is our domain, not the domain of the computer.

%-- not to experiment with examples, or to do a lengthy case check, but to actually \emph{prove or check proofs of theorems}. In particular, it is a situation where computers can actually wrestle with the infinite. Furthermore this ``new'' technique applies \emph{across pure mathematics}, not just to areas in which computers have been traditionally used in the past.
However, this is no longer the case. Computer proof assistants are programs which know the axioms of mathematics. A consequence of this is that they can do both computing in the traditional sense, and also \emph{reasoning}. In practice this means that one can write some computer code in a proof assistant which corresponds to the proof that there are infinitely many primes~\href{https://leanprover-community.github.io/lean-web-editor/#url=https%3A%2F%2Fraw.githubusercontent.com%2Fkbuzzard%2Fxena%2Fmaster%2Fsrc%2FICM%2Finfinitude_primes.lean}{(\underline{live link})},
or even to a proof~\cite{dhl} of the main result in a recent Annals paper~\cite{ellenberg-gijswijt}.

I wrote ``new'' in quotes above because it is not new at all; computer scientists have been creating tools like this for decades now. Indeed, the first computer proof assistants appeared in the 1960s. However, more recently three things have happened. First, the technology has now reached the point where research level results across all of the traditional mainstream areas of pure mathematics are now simultaneously accessible to these systems, at least in theory, and, increasingly, in practice. Secondly, the systems are far more autonomous than they used to be. Tactics are commands which can be designed by users and which are capable of putting together hundreds if not thousands of tedious axiomatic steps, enabling mathematicians to communicate with these machines in a high-level way, similar to the way which they communicate with each other. Finally, and crucially, research level mathematicians are finally beginning to get involved; we are seeing material at MSc level and beyond being formalised, by mathematicians, across many areas of mathematics now. These developments mean that teaching research level material to a computer proof system in all areas of mathematics is now becoming a feasible possibility -- indeed, it is already happening right now, and shows no signs of stopping.

%The foundations of the subject lie on the boundary between two fields -- pure mathematics, and theoretical computer science. Over the last few years there has been a large amount of synergy and cross-fertilization. Computer scientists beginning to understand the kind of objects modern mathematicians use \emph{across pure mathematics}, and conversely mathematicians are beginning to learn how best to embed the objects they know and love such as schemes and manifolds within a computer proof system. Of course for simpler objects such as groups and rings, this kind of thing was understood decades ago; what is happening now is that because more and more young mathematicians are engaging with the software, results involving far more complex objects are being formalised, and conversely computer scientists are beginning to get a much better feeling for the kind of things which a modern mathematician needs to be able to do easily, and are designing their theorem provers accordingly.

%Of course, expecting a computer to come up with a proof of the Riemann Hypothesis any time soon would be absurd. In fact the reason the chess/go analogy is a weak one is that even though chess and go have gigantic game trees, they are \emph{finite}. Mathematics takes off the moment you lay down the axiom or construction which enables the natural numbers -- an infinite object -- to be constructed, and from this point on we are playing a different game, and one which computer scientists know much less about.

This paper consists of 4 sections, which are independent of one another, and can be read in any order.

The first is historical; it consists of descriptions of the systems which are being, or have been used, to formalise mathematics, and discussions of results which have been taught by humans to computers over the last 20 years. It also notes various historical technical advances.

The second is an overview of one of the largest currently available monolithic mathematical libraries in existence, namely Lean's mathematics library {\tt mathlib}. Lean~\cite{leanpaper} is a free and open source computer proof assistant written primarily by Leonardo de Moura at Microsoft Research. Lean's maths library {\tt mathlib}~\cite{mathlibpaper} is a free and open source library for Lean, developed by a community of users across the world, ranging from undergraduates to professional mathematicians. {\tt mathlib} is the library which has powered several of the most recent significant results in the area.

The third section consists of an introduction to type theory as a foundation for mathematics; it explains how mathematical structures, theorems and proofs can be encoded within these foundations. Note that many of the modern computer proof systems where non-foundational mathematics is happening (Lean, Coq, Isabelle/HOL) use type theory rather than set theory; however type theory proves the same theorems as set theory. Furthermore, mathematicians who can prove theorems but who do not know the axioms of ZFC set theory can happily write code in a type theory proof system corresponding to these theorems without knowing the axioms of type theory either.

Finally, a speculative final section describes in more detail some personal ideas of the author and others about the kinds of things which software such as this can be used for, and how it might help us to do our jobs.

\subsection*{Acknowledgements}

I thank the Lean prover community for welcoming me, a mathematician with very little programming experience, into their community back in 2017, and also for reading and giving extensive comments on a preliminary version of this article. Patrick Massot in particular sent many helpful comments on a first draft. I thank Assia Mahboubi and Manuel Eberl for giving advice on the Coq and Isabelle/HOL code in this paper, and to both them and Jeremy Avigad for helpful historical comments. Finally I would like to effusively thank Leonardo de Moura for writing my favourite computer game, and Mario Carneiro for teaching me how to play it.

\tableofcontents

%One of the things I describe in this paper are some more reasonable goals for ``computer-powered mathematics'', some of which are already attainable. The goals are releated to both teaching and research. Each goal represents genuine progress, in the sense that it will be a new a way in which computer proof systems can help the modern mathematician to do their job. Whether computers will be replacing human mathematicians any time soon is a question for the science fiction writers. But how much they will be helping us in the near future depends (to borrow a phrase from Jeremy Avigad~\cite{avigadskiningame}) on how much skin we are prepared to put in the game, as a community.

%Here is an overview of what I will discuss.

%1) history of 21st century formalisation attempts. Mizar, Isabelle/HOL, Coq,

%2) Example of code

%3) snapshot of mathlib. Lean's maths library learning at approx the same rate as an undergraduate.

%4) LTE

%5) Future goals. where to go next and why? [stacks]

%Appendix) sets, types, topological spaces up to homotopy: very good for homotopy theory

\section{A brief history of formally verified theorems}

%Let me finish this introductory section by explaining why formalisation of mathematics is harder than you might think.

In this section I will talk about the previous successes of computer proof assistants -- computer programs which check human proofs -- in mathematics. There are far more projects here which I could have mentioned, and I apologise to those who have undertaken major mathematical formalisation projects which I have not cited. Examples of computer proof assistants in which a substantial amount of mathematics has been formalised include Lean~\cite{leanpaper}, Coq~\cite{coqpaper}, Isabelle/HOL~\cite{isabellepaper}, HOL Light~\cite{hollightpaper}, Metamath~\cite{metamathpaper} and Mizar~\cite{mizarpaper}.

For a computer to formally verify a theorem, it ultimately needs to be able to deduce the theorem from the axioms of the foundational system (typically set theory or type theory) which the proof assistant has been designed to use. I will use the below discussion of historical results to introduce some conceptual breakthroughs which have over the years enabled the formalisation of mathematics to become feasible.

This section cannot do justice to all of the work which has occurred in the area; I thoroughly recommend Hales' paper ``Mathematics in the age of the Turing machine''~\cite{hales-turing} for more background and examples, although much has happened since that paper was written in 2014.

\subsection{The 20th century}

Consider the problem of proving from first principles that if $x$ and $y$ are real numbers, then $(x+y)(x+2y)(x+3y)=x^3 + 6x^2y + 11xy^2 + 6y^3$. We all know that the real numbers are a commutative ring, so let us assume that fact. The question now becomes how to use the axioms of a commutative ring to prove the equality that we want. How many lines would a proof from first principles be? Surely not too many! We apply distributivity a few times to expand out the brackets on the left hand side, and then of course it just becomes a matter of tidying up and equating terms. As humans we do not think too much about the tidying-up process, however if you try proving this in a theorem prover then you will discover that actually it is a combinatorial nightmare. For example there is a step in the proof where where we need to prove something of the form
$$((A+B)+(C+E))+((D+F)+(G+H))=((((((A+B)+C)+D)+E)+F)+G)+H$$
using only the laws of commutativity and associativity of addition. Humans apply a \emph{principle} to justify this step, not an axiom, and indeed proving such a triviality using only the axioms of a ring is surprisingly fiddly. There is also the issue of turning things like $x((2y)x)$ into $(2(x^2))y$ and so on. 

The very early theorem provers had very limited ability to apply principles, meaning that proving results such $(x+y)(x+2y)(x+3y)=x^3 + 6x^2y + 11xy^2 + 6y^3$ would need to be done manually, meaning something like a 30 line proof. If such a triviality hides 30 lines of axiomatic mathematics, imagine what is hidden behind claims of the form ``The function $f$ is clearly $O(x^{-2})$ for $x$ large''? It is one thing writing a computer proof assistant -- it is quite another one to write one which scales to do the kind of things which we humans do intuitively. For this and other reasons, many of the earlier formalisation achievements of the 20th century were mathematically trivial. In particular, there were many proofs of the irrationality of $\sqrt{2}$ and of the infinitude of primes, but these were being used as benchmarks for the systems.

In the final two decades of the 20th century, computer provers began to appear which had new functionality. In these later systems, users could write ``tactics''. Tactics are computer code which assembles axiom applications together into principles. For example in a modern prover like Lean, $(x+y)(x+2y)(x+3y)=x^3 + 6x^2y + 11xy^2 + 6y^3$ can now be proved in one line by invoking the {\tt ring} tactic\footnote{See~\cite{assia-ring} for a description of the sort of issues which arise when writing such a tactic.}. Tactics allow formalised mathematics to more closely resemble ordinary mathematical practice by making "obvious" things automatic.

\subsection{The prime number theorem}

In 2004, a team comprising of Jeremy Avigad, Kevin Donnelly, David Gray and Paul Raff formally verified the prime number theorem, in the Isabelle/HOL system. The proof they formalised was the Erd\H{o}s--Selberg ``elementary proof''. The work used inputs from both arithmetic and basic real analysis. Of course calculations involving growth of functions which look easy on paper still took time and effort to formalise. Manipulation of inequalities which to humans look easy need to be done either by hand or via a Fourier--Motzkin elimination tactic in a theorem prover. The reason that the Erd\H{o}s--Selberg proof was preferred to the traditional complex analysis proof was that at that time Isabelle/HOL had no complex analysis library at all. 

What we conclude is that by 2004, more serious undergraduate and MSc level material was now in theory accessible to these systems, at least in some areas of mathematics. We also see that we are at a stage where libaries of proofs in distinct areas of mathematics are able to interact with one another.

In 2009 John Harrison formalised the complex-analytic proof of the prime number theorem in the HOL Light theorem prover~\cite{harrison-pnt}, motivated in part by the fact that HOL Light already had a theory of complex analysis including Cauchy's integral formula. In 2016 Mario Carneiro formalised the Erd\H{o}s--Selberg proof in Metamath, a set theory based prover which has essentially no tactics; as you can imagine this was a heroic effort.

Thus the Prime Number Theorem became some kind of a poster child for formalisation. One can understand why -- it was a celebrated theorem in mathematics, the proof is not at all trivial, and any formalisation in a theorem prover demonstrates that the prover is capable of reasoning about both the discrete and the continuous simultaneously.

As may be becoming apparent to the reader however, one reason that the result was being independently formalised in several theorem provers was that it is extremely difficult to translate a proof written in one of the systems to a proof in another system. One issue is that different systems might have different foundations; for example HOL Light and Isabelle/HOL are type theory systems, and Metamath is a set theory system. Another issue is that even if two proof systems have very similar foundations, they might have different \emph{idioms}; different libraries in different systems could be set up to do the same thing in very different ways. Without getting too technical, in order for these computer proof systems to work one has to have some kind of a method for moving between structures ``behind the scenes'' -- for example the reals are a field, and hence they are an additive group (and a multiplicative monoid), and in particular one wants all theorems about additive groups such as $0+a=a$ to apply instantly to fields such as the reals without any fuss. Humans of course have no problems with this, but in a computer proof system one needs some kind of infrastructure which is making this happen automatically, and if different systems are doing this in different ways then of course this makes automatic proof translation much harder.

Thus it came as a shock to me when in 2020 Mario Carneiro announced that he had used his Metamath Zero project~\cite{mm0} to port the Metamath proof of the prime number theorem to Lean. The two systems are about as far apart as it is possible to be -- Metamath uses set theory as a foundation and Lean uses type theory, for example. Metamath proofs are typically far more low-level, with limited automation available, whereas typical Lean proofs are very tactic-heavy. However the system worked, and produced code which compiled; it was of course also unreadable. It was tens of thousands of lines of completely unmotivated primitive code defining variables and applying basic principles of logic, with no comments. In fact it was a wonderful example of something which satisfied a formal definition of ``being a proof'', whilst in some sense imparting no information whatsoever to the human reader other than the fact that the theorem was true.

Of course if computers begin to write proofs by themselves, they might all look like this, at least at first.

\subsection{The four colour theorem}

The four colour theorem (formerly the four colour conjecture) was a notorious problem in graph theory raised in the 1850s and which remained unresolved for over 100 years. One formulation of it is the assertion that the vertices of every planar graph can be coloured with four colours in such a way that no two adjacent vertices share a colour. The statement is an elegant combinatorial problem, and it came as a shock to some in the mathematical community that the proof, announced by Appel and Haken in 1976, used a computer in an essential way. Appel and Haken constructed a collection of 1834 graphs with the property that a minimal counterexample must contain one of these graphs as a subgraph, but that conversely no graph containing one of these 1834 graphs as a subgraph can be a minimal counterexample. The verification of these claims was done using a bespoke computer program which, in those days, took over a month to finish running. The Appel--Haken proof was an outlier because whilst the principle of the proof was possible to understand, the details were too difficult for a human to follow in practice; one billion case splits (this is what the computer part of the proof looks like) is not something which humans can do manually and accurately within a reasonable time frame. The proof relies, essentially, on a computer calculation and hence it relies on the correctness of the computer code. Small bugs in computer code are of course commonplace, although it of course could be added that small bugs in human-written proofs are also commonplace. However the mathematical community is well-equipped to discover and fix small bugs in human-written proofs, and was perhaps rather less well-equipped to verify the correctness of computer code, especially in 1976.

In 2004 Georges Gonthier finished a formal verification of the Appel--Haken result -- more precisely he formalised the 1997 Robertson--Sanders--Seymour--Thomas variant of the argument~\cite{gonthier-4ct}. The work comprised of 60,000 lines of code written in the Coq proof assistant. In particular it completely dwarfs the prime number project discussed in the previous section. It contains a complete formalisation of the theoretical part of the work -- formal proofs of results in topology (to reduce the statement about arbitrary planar graphs to one of a discrete combinatorial nature) and graph theory -- whilst also formally verifying the computer calculation necessary to finish the proof. Note in particular that (as in the proof of the prime number theorem) much of the work comprised of writing foundational material rather than formalising the proof itself.

It is interesting to note that the process of formalisation led to simplifications in the argument. For example Gonthier developed a theory of what he called combinatorial hypermaps, which greatly reduced the amount of topology needed in the proof, and in particular removed the dependency of the argument on the Jordan Curve theorem. Gonthier developed some original mathematics as part of the work -- for example he isolated a combinatorial criterion for his hypermaps which was equivalent to planarity.

Naively, it looks like in this case we are replacing one ``proof by computer'' with another one, however this is missing the point. Firstly, the Coq formalisation covers not just the Appel--Haken computer code, but also all of the rest of the Appel--Haken argument. Secondly, one can view the formal verification as an independent check of the proof. Finally, instead of having to trust the code written by Appel and Haken and which few people have read, we are instead having to trust code written by the authors of Coq. Coq has been around for a long time (the first version was written in 1984), has a small kernel, and the system has many users. A bug which meant that Coq could incorrectly claim that an unproved theorem was true would be unlikely to manifest itself in just one project and is far more likely to ultimately be discovered. In contrast, the Appel--Haken code is a bespoke piece of code with few users so arguably bugs are more likely.

Gonthier wrote a very informative piece~\cite{gonthier-notices} about his work for the Notices of the American Mathematical Society (including an exposition of the theory of hypermaps), as part of the November 2008 issue; this issue was devoted to formal verification of mathematics in a computer proof system and provides an excellent survey of the field as it then stood.

\subsection{The odd order theorem}

The odd order theorem is the theorem that any finite group of odd order is solvable. In 2013 a team of 15 people led by Gonthier formally verified a proof of this theorem in Coq~\cite{gonthier-oddorder}. This piece of work is notable for several reasons. Firstly, the proof is very long; a complete argument (modulo the basics in group and representation theory) is presented in the two volumes~\cite{bender-glauberman} and~\cite{peterfalvi}. Secondly, we have moved way beyond MSc level mathematics here -- this work was one of the reasons that Thompson was awarded the Fields Medal in 1970. The proof is a very delicate argument in finite group theory, much of which involves analysing the structure of a minimal counterexample and ultimately showing that it cannot exist. The Coq proof involved formalising both of the books mentioned above, plus of course all the background material in group theory, representation theory, Galois theory and number theory; indeed, formalisation of the background material took up much of the six years which the authors spent on the proof. Figuring out how to handle such a large-scale formalisation project was also a non-trivial task.

It is perhaps worth stepping back and asking how work like this contributes to human understanding. The naive answer to this is ``it guarantees that the human proof is correct''. However, in my opinion this is not the main contribution. Humans were well aware even back in the 1960s that the proof was correct -- had there been any doubt, Thompson would not have got the Fields Medal. What the formalisation work shows us that theorem provers have now become able to operate at this kind of scale. Entire books of mathematics can now be formalised in one system without the system running out of memory or grinding to a halt. On average, one line of mathematics in~\cite{bender-glauberman} or~\cite{peterfalvi} corresponded to five lines of computer code, so we learn that by 2013 the so-called ``De Bruijn factor'' for this kind of mathematics is around~5. However this ratio should not be taken too seriously: in parts of the argument the ratio is essentially one, and in other parts it is much larger. Note also that this factor may vary considerably between theorem provers.

We also learn that large formalisation projects such as this are a very effective way to motivate development of foundational mathematics libraries. One consequence of this formalisation project was that Coq developed a very solid library of undergraduate-level algebra which can of course be used (and is used) for other projects.

The write-up~\cite{gonthier-oddorder} of the odd order work is an interesting read. Some sections concentrate on the mathematics or the history, but there is also a discussion about constructive mathematics, something which I felt should have nothing to do with the work, and also about implementation issues, something else which mathematicians typically do not ever have to think about. For example, one observation made in section 3 was that many theorems involving two or more finite groups would usually be formalised assuming that these groups were both subgroups of some larger ambient finite group~$X$. This can be done without any loss of generality of course, because given two groups~$G$ and~$H$, they are both subgroups of $G\times H$. Why is this observation important? This is an \emph{implementation issue} -- the domain of the computer scientist. Working with subgroups rather than groups might be easier, or nicer, when it comes to actually implement certain theorems in the theorem prover. It is worth noting however that such a trick does not work more generally: for example in algebraic geometry one uses the category of commutative rings with 1, and morphisms by definition send~1 to~1. If~$R$ and~$S$ are general commutative rings with~1 then there is in general no morphism of rings from $R$ to $R\times S$ sending~1 to~1, so one is forced to implement commutative ring theory in a more ``traditional'' manner. See~\cite{eric-scalar} for how this was done in Lean.

Regarding constructivism -- the authors of the work put in a lot of effort to keep their proof ``constructive'', for example the avoidance of all uses of the complex numbers when setting up the basics of representation theory. The complex numbers do not have decidable equality, meaning that there is in general no algorithm for proving that two constructively defined complex numbers are equal (for example, one can evaluate a definite integral numerically and observe that it seems to be $0$ to 1000 decimal places, but there is not some generic algorithm which we can apply to an arbitrary integral in order to decide whether or not it equals $0$). This means that in constructive mathematics, where the law of the excluded middle cannot be assumed, one cannot do a case split on whether $z=0$ or not, if $z$ is a complex number, and more generally plenty of constructions become noncomputable and hence much harder to reason with constructively. These design choices thus increase the amount of work needed to get representation theory working. I had thought that constructivists had died out in the early part of the 20th century. It turns out that they are alive and well, and typically working nowadays in computer science departments. One reason for this is that constructivism plays an important role in the theory of programming languages. Reluctance to use the law of the excluded middle is to a certain extent a cultural decision. However there are also situations where working constructively enables a computer proof system to prove certain results ``automatically'' (for example by an explicit computation). Whilst working constructively may have been feasible for a project about finite groups, the law of the excluded middle is used throughout most modern research level mathematics and it is not really feasible to work constructively when doing the kind of mathematics which is happening nowadays in mathematics departments. However it is also worth stressing that most modern proof assistants have no problems with the law of the excluded middle, the axiom of choice, and other non-constructive axioms -- they are available, if you want them. Certain proof strategies are ruled out if one chooses to work nonconstructively, but one can counter this by writing new tactics specifically designed to do computations in fields such as the real and complex numbers. Nonconstructive axioms are used extensively in Lean's mathematics library {\tt mathlib}, for example.

\subsection{The Kepler conjecture}

The Kepler conjecture states that the face-centred cubic packing is the densest way to pack congruent spheres in 3-space. Hales and Ferguson proved the conjecture in 1998; it had at that point been open for over 350 years (it was raised by Kepler before Fermat proposed his Last Theorem). Part of the Hales--Ferguson proof involved the checking of over 23000 non-linear inequalities on a computer; another part involved a computer classification of all tame graphs. Other computer calculations were also involved. In this respect the proof is similar to the Appel--Haken proof of the four colour theorem; computations need to be carried out which are simply far too great for humans to do in a reasonable time frame.

Because the result was regarded as important, the referees felt duty-bound to attempt to check the computer part of the proof in some way; however ultimately they gave up, and in~\cite{hales-kepler} Hales states that the paper was published (in the Annals) without complete certification from the referees. In 2003 Hales announced a project to formally verify the proof using computer proof systems. Hales used a combination of HOL Light and Isabelle/HOL, and the project turned into an international collaboration, with 22 authors listed on the final paper. The formalisation project took around 12 years to complete, and comprised over half a million lines of code. Just as for the other projects in this section, one of the main benefits of the work to the formal proof community is that HOL Light's standard library grew to include theorems such as the Brouwer fixed-point theorem, the Krein--Milman theorem and the Stone--Weierstrass theorem.

In 2017 Hales gave a talk~\cite{hales-newt} at the Newton Institute where he tells the story of the Kepler proof, and explains a vision for the future of formalised mathematics. This talk was, for me, the turning point, and was one of the main motivations behind the work described in the following subsection.

\subsection{Perfectoid spaces}

The previous formalised results all have something in common. Whilst some of them represent truly deep mathematics, all of the formalised proofs involve reasoning about objects which are in some sense \emph{elementary} (planar graphs, prime numbers, finite groups, spheres). Furthermore, most (but not all) of the formalising done prior to 2017 was being done by computer scientists. In Hales' talk linked to above, he coherently argued that for further progress in this area, this state of affairs had to change. At that time I had only just begun to dabble with computer proof assistants and my initial plans were to attempt to integrate them into my undergraduate teaching. However Hales' arguments resonated with me, and within a few months I found myself working with undergraduates at Imperial College, formalising the definition and basic properties of schemes in the Lean theorem prover. This project involved developing basic theories of localisation of rings and of sheaves on topological spaces; however it was relatively straightforward (modulo poor design decisions; the reader interested in more details can see them in~\cite{ic-schemes}). I was thus shocked to discover afterwards that schemes -- such a basic notion in algebraic geometry -- had not been previously formalised in \emph{any} other computer proof system! Furthermore, the project made it quite clear to me that formalising far more heavyweight mathematical objects should easily be possible.

In late 2017 Patrick Massot (a topologist) and myself independently came up with the idea of formalising perfectoid spaces; the topic was in the air because it was at that time an open secret that Scholze was going to be awarded a Fields Medal for his invention/discovery of the concept and its applications to arithmetic geometry. I knew the mathematical definition, having dabbled in the area myself, and when Johan Commelin, another arithmetic geometer, appeared in the Lean Zulip chatroom in 2018 the three of us decided to go for it. Around 16000 lines of code and eight months later, we had a formalised definition; one could summarise the work as a computer formalisation of the single line of mathematics ``let $X$ be a perfectoid space''.\footnote{In the odd order formalisation, the de Bruijn factor (ratio of lines of computer code to lines of human text) was around~5. Here one could argue that it is 16000. However one could of course also argue that it might well take several thousand lines of human text to define a perfectoid space in full.}

The work was of course partly intended as a public relations stunt; computer scientists were well aware of the existence of computer theorem provers, however mathematicians seemed not to be, and this was an attempt to make them notice. The plan was a success -- the project did seem to raise the profile of computer theorem provers within the mathematics community. Note however that we did not construct any examples of perfectoid spaces other than the empty perfectoid space\footnote{To prove that the empty set can be given the structure of a perfectoid space, one needs to check that an arbitrary product of trivial topological rings is the trivial topological ring.}, and all three of us were well aware of the problems preventing us from formalising any of Scholze's serious theorems about perfectoid spaces at that time; we were missing so many of the prerequisites. As with previous projects, one tangible gain from the work was the growth of the mathematics library of the system in question. Most of the results in Bourbaki's General Topology ended up as part of Lean's mathematics library {\tt mathlib} as a result of this project, as well as plenty of results in topological algebra, and it also motivated the beginnings of a theory of valuations and discrete valuation fields.

One can consider the perfectoid space work as in some sense being orthogonal to what was usually being attempted in a theorem prover. Many of the prior results highlighted in this section are proofs of long and complex theorems about relatively simple objects. The proof that the empty set can be given the structure of a perfectoid space is a very simple theorem about a much more complex concept. Of course the natural next question is whether computer proof systems can prove complex theorems about complex objects. One year after the perfectoid space project, we began to find out.

\subsection{Condensed mathematics}

Clausen and Scholze have been developing a theory of condensed mathematics. A condensed set is a variant of a topological space. The main insight is that condensed objects may have better homological properties than topological objects (for example the category of condensed abelian groups is an abelian category, whereas the category of topological abelian groups is not). They hope that these ideas will enable techniques in homological algebra to apply to new areas of analytic geometry. At the end of 2020, Scholze approached me and asked if we had had a study group on the work at Imperial; I answered that we had. Scholze then asked whether we had looked through all the details of the proof of Theorem~9.1 of~\cite{scholze-analytic}; I answered that we had not. Scholze then remarked that he had had the same response from other mathematicians, and raised the possibility that perhaps nobody other than himself and Clausen had ever read the proof carefully. Furthermore he suggested that perhaps this might remain true even after the refereeing process. The reason he was concerned about this was that, for Scholze, this was the theorem that the entire theory stood upon. The proof was very technical; it built upon a more ``elementary'' but rather unwieldy intermediate result, Theorem~9.4 of~\cite{scholze-analytic}. Scholze agreed to challenge the formalisation community to prove his Theorem~9.1 in a blog post~\cite{xena-scholze1}, later published as~\cite{scholze-expmath}. Although the challenge was to the formalisation community in general, it seems that only the Lean community responded; this is perhaps unsurprising, as (for perhaps only for sociological reasons) it has come to be the case that mathematicians interested in ``the kind of mathematics which wins Fields Medals'' and also interested in theorem provers tend to gravitate towards Lean.

Johan Commelin became the de facto leader of the formalisation process, with Patrick Massot supporting him in making a blueprint~\cite{lte-blueprint} of the strategy (that is, a carefully-written roadmap) and a team of algebraic number theorists, arithmetic geometers and other mathematicians (Riccardo Brasca, Damiano Testa, Filippo Nuccio, Adam Topaz, myself, Patrick Massot, Bhavik Mehta\ldots) then began working on the project, with the occasional help from people with a computer science background such as Mario Carneiro. Within six months the team had grown to over ten people and we had formalised a complete proof of Theorem~9.4 (see~\cite{nature}). At the time of writing we have not deduced Theorem~9.1, but it is only a matter of time. A second blog post~\cite{xena-scholze2} by Scholze indicates his thoughts on the matter; in particular we see that he is now far less concerned about the situation regarding the correctness of the results. Furthermore, Scholze has indicated (personal communication) that the process has enabled him to better understand what powers the proof, and Commelin not only learnt the mathematics as part of the process, but also simplified the argument in several places, most notably in the removal of the dependency of the argument on prior work of Breen and Deligne.

For me, this represents substantial evidence that now \emph{any pure mathematics} can be formalised in theorem provers -- both in theory, and in practice. It takes time, but it is possible. The formalisation of the work led both to better understanding of it, and to simplifications of the argument. Also worth mentioning is that, as in many other formalisation projects, a substantial amount of time was spent formalising background material (for example the theory of normed groups and the theory of profinite spaces). As the libraries of the provers get better and start to contain the kind of material which working mathematicians take for granted, there will be fewer of these ``startup costs''.

\subsection{Other results}

There are plenty more examples of serious formalisation efforts which we do not have the space to cover. We list some examples here. Gou\"ezel formalised the basic definitions of $C^k$ and $C^\infty$ manifolds in Lean, extending earlier work done in Isabelle/HOL. Mahboubi and Sibut-Pinote proved irrationality of $\zeta(3)$ in Coq~\cite{mahboubi-zeta} and Eberl proved it in Isabelle/HOL~\cite{eberl-zeta}. Mahboubi has also done extensive work on rigorous numerical values of integrals in Coq, and also in Coq Bertot, Rideau and Th\'ery formally verified the first one million decimal digits of~$\pi$~\cite{pi}. Eberl has formalised much of Apostol's textbook on analytic number theory in Isabelle/HOL~\cite{eberl-apostol}. Han and van Doorn proved independence of the continuum hypothesis in Lean~\cite{flypitch}. Immler formally verified Tucker's calculations used to verify the existence of the strange attractor~\cite{immler-strange}. Mehta and Dillies formally verified Szemer\'edi's regularity lemma and Roth's theorem on arithmetic progressions in Lean, and Edmonds, Koutsoukou-Argyraki and Paulson verified them in Isabelle/HOL~\cite{isabelleRoth}. The Poincar\'e--Bendixson theorem was formalised by Immler and Tan~\cite{poincarebendixon} in Isabelle/HOL; note that the usual proof as understood by mathematicians relies on drawings, and formalising drawings can be hard work. The Ellenberg--Gijswijt resolution of the cap set conjecture was verified in Lean by Dahmen, H\"olzl and Lewis~\cite{dhl}. Commelin and Lewis constructed Witt vectors and showed that $W(\mathbb{F}_p)=\mathbb{Z}_p$ in~\cite{witt}; this work is interesting because not only did they formalise the delicate mathematics involved, they also wrote tactics which would enable them to reduce various questions to the universal case in a painless manner. Finiteness of the class group of a global field was proved in Lean by Baanen, Dahmen, Narayanan and Nuccio in~\cite{class-group} (it still astonishes me that this result, special cases of which were known to Gauss and which is a standard theorem in an undergraduate mathematics degree, was formalised in a proof assistant for the first time in 2021; this development is a demonstration of how the mathematical interests of the formalisation community now is encompassing work which people in mathematics departments might view as ``mainstream'').

There is also work in progress (at the time of writing). Teams of people who collaborate on the Lean Zulip chat~\cite{zulip} are currently working on a proof of Fermat's Last Theorem for regular primes, and on Smale's theorem that it is possible to evert a sphere. A general project to formalise many basic results in the theory of schemes is also underway.

\section{{\tt mathlib}}

In this section I will give an overview of Lean's mathematics library, one of the largest monolithic collections of formalised mathematics in existence and, more importantly, one which is currently experiencing rapid growth. To a certain extent it is a personal perspective; a different point of view, which talks more about the computer science powering the library, is presented in~\cite{mathlibpaper}. 

The principal developer of the Lean Theorem Prover is Leonardo de Moura, who started the project in 2013. At the 2017 Big Proof conference in Cambridge it was decided to split off most of the ``mathematical'' part of the prover from the ``core'' part, and move the mathematics into a library of its own. Thus {\tt mathlib} was born. At the time {\tt mathlib} contained definitions of groups, rings and topological spaces, filters, a construction of the rational numbers (the naturals and integers remained in core Lean), and little else. Johannes H\"olzl and Mario Carneiro became the maintainers of the library, and between them they began to slowly build more mathematics, for example, the real numbers. H\"olzl had written a lot of the topology part of the repository, following the Isabelle/HOL approach which relied heavily on the concept of a filter. Carneiro wrote a robust theory of finiteness, and slowly the library began to become relevant to the ``working mathematician''.

The library is a free and open source project. It is monolithic in the sense that there is one definition of a group, one definition of a ring, one definition of the real numbers and so on, and all of these definitions can be imported simultaneously and interact with each other. Initially it was not clear what its goals were, other than being a place where people could experiment with doing mathematics in Lean. Mathematicians such as Scott Morrison, myself and Patrick Massot got involved at a very early stage, and because our background was in mathematics which relied on classical logic (i.e., the law of the excluded middle) and other non-constructive axioms such as the axiom of choice, the library developed with these classical assumptions at its core. Each successful mathematics project written in Lean and powered by {\tt mathlib} seemed to attract more mathematicians to its chatroom, which in turn led to more projects. Within a couple of years Lewis had formalised the $p$-adic numbers~\cite{lewis-padic}, myself and a team of undergraduates (Lau, Hughes, Livingston and Fern\'andez Mir) formalised schemes~\cite{ic-schemes}, Dahmen, H\"olzl and Lewis formalised the 2017 Ellenberg--Gijswijt Annals proof of the cap set conjecture~\cite{dhl}, and Massot, Commelin and myself formalised the definition of a perfectoid space~\cite{bcmperfectoid}. Each of these projects could not have happened without {\tt mathlib}; conversely each of these projects contributed to the growth of {\tt mathlib}.

Plenty of developments were also taking place which were not written up as papers, and whose main purpose was simply to grow {\tt mathlib}. I supervised student projects where undergraduates could formalise material they were learning in class and add it to the library; for example Sylow's theorems (Chris Hughes), nilpotent groups (Ines Wright), conformal maps (Yourong Zang) and the Radon--Nikodym theorem (Kexing Ying) were added this way. Amelia Livingston developed a theory of localisation of monoids and rings which we needed for algebraic geometry. I pushed undergraduates (Hughes, Lau, Lee) to formalise a standard Galois theory course in Lean; they developed a theory of field extensions, and the project was then taken up by a group of mathematics PhD students in Berkeley (Miller, Browning, Lutz) who finished the job, proving the fundamental theorem of Galois theory and the insolvability of the quintic~\cite{browninglutz} (note that this was coincidentally formalised in Coq just a couple of months beforehand~\cite{assia-quintic}). Baanen, Dahmen, Narayanan and Nuccio formalised a proof of the finiteness of the class group of a global field~\cite{class-group}. I was pushing algebra, but others were pushing geometry and analysis. Gou\"ezel and Macbeth developed a theory of manifolds, and Gou\"ezel and Kudryashov developed an extensive theory of single and multivariable calculus, including the implicit function theorem and the Picard--Lindel\"of theorem. Gou\"ezel also formalised the Gromov-Hausdorff space: a metric space parametrising nonempty compact Hausdorff metric spaces up to isometry.

Morrison has developed a huge amount of category theory, and he and Topaz have now formalised the definitions of abelian categories and the beginning of the development of derived functors and homological algebra. Massot has developed valuation theory and a theory of completions of uniform spaces and of topological groups and rings. Tuma developed the theory of Jacobson rings, and I developed some of the basics of other standard ideas in commutative algebra (projective and flat modules, discrete valuation rings), and Springer, Kuelshammer, and many others have also contributed to algebra. H\"olzl developed the theory of Lebesgue measure, and van Doorn formalised Haar measure. There are many more people who have made contributions ({\tt mathlib} now has over 200 contributors) and new contributions are always welcome. Contributions are reviewed by the maintainers. One of the principles of the library is to do things ``in the correct generality''. This meant, for example, that multivariable calculus and some exotic integrals taking values in Banach vector spaces was developed first, and single variable calculus was deduced as a corollary. The library is not optimised for pedagogy or readability; the idea is to continue to make a solid foundation for the kind of mathematics which is happening in a contemporary mathematics department. 

It is interesting to note that Lean seems to be learning mathematics at around the same speed as an undergraduate. In the four years which the library has been growing, it has gone from essentially zero to a solid MSc level coverage in number theory and commutative algebra, and BSc level real analysis. In complex analysis, differential geometry and representation theory it is perhaps not quite yet at final year BSc level, but things move fast and this sentence, written in 2021, will quickly date. For an up to date idea of the current status of {\tt mathlib} the best idea is to take a look at the Lean community's full overview of mathlib~\cite{community-mathliboverview}, or its summary of the undergraduate level mathematics it contains~\cite{community-ugmaths}.

\section{A brief guide to type theory}

In this section we explain the basics of type theory and how it can be used as a foundation of mathematics. Many modern theorem provers use some version of type theory as their foundations. For example Isabelle/HOL and the other HOL systems use simple type theory, Lean and Coq use dependent type theory, and the various HoTT systems developed by Voevodsky and others use homotopy type theory. There are a few computer proof assistants which use set theory -- Metamath and Mizar are the two most prominent -- however it is not unfair to say that nowadays most mathematics done in theorem provers is done in a type theory system, so a mathematician interested in dabbling in formal proofs should at least know something about the basics, which is what this section attempts to describe.

\subsection{What is a type?}

Mathematicians nowadays are used to seeing the word ``set'' floating around when it comes to basic definitions. For example, we are told that a group is a set equipped with a multiplication such that some axioms hold. We're not told what a set is though; a course on ZFC set theory tells us a list of properties which sets \emph{have}, but they don't tell us what a set \emph{is}. Indeed in this context the word ``set'' has no formal definition; it is simply the generic term for an object in our model of the axioms of mathematics, and we build other mathematical objects on top of this basic object.

In definitions such as the definition of a group, the word ``set'' is being used to mean no more than ``collection of elements''. In type theory, the role of a ``collection of elements'' is played by the \emph{type}. A type is a collection of terms. The definition of a group in type theory: a group is a type equipped with a multiplication such that some axioms hold. The only difference is the notation: the set-theoretic $a\in X$ is replaced by the type-theoretic \lstinline{a : X}.

As mentioned above, those of us who have been to a set theory class will know that, when using set theory as a foundation of mathematics, \emph{everything} is a set. For example the elements of a group are, strictly speaking, also sets, so one could in theory talk about their elements too, although within the context of group theory such questions would not be mathematically meaningful, as they are not isomorphism-invariant. In type theory this is not possible; the elements of a type are called \emph{terms}, and in general terms are not types. In type theory, everything is a term, and every term \emph{has} a type, but not every term \emph{is} a type. For example, in type theory $37\pi^2$ is a term, whose type is $\mathbb{R}$, the type of real numbers. We write $37\pi^2\,:\,\mathbb{R}$. However $x\,:\,37\pi^2$ \emph{does not make sense}, because $37\pi^2$ is not a type. In a set-theory based theorem prover, questions such as asking if the trivial group is an element of the Riemann zeta function would make sense but its meaning would be unmathematical -- it would depend on implementation decisions. Type theory thus provides a basic check that what you are writing has mathematical meaning.

In a type theory system, the type $\mathbb{R}$ is still built from $\mathbb{Q}$ as equivalence classes of Cauchy sequences, or via Dedekind cuts, or as another of the standard constructions; the mathematical part of the story is identical to the set theory set-up, it is just that the language used is slightly different (types and terms, rather than sets and elements).

One difference between types and sets however is that \emph{types don't mix}: distinct types are disjoint. This has practical advantages when formalising mathematics because it provides a strong check that the mathematics you're typing \emph{makes sense}: in type theory, if $g$ is an element of the group $G$, then the only type that $g$ can ever be a term of is~$G$.

This approach does however have consequences which can initially come as a shock to a mathematician. For example one could make a type representing the positive reals $\mathbb{R}_{>0}$ and a type representing the reals $\mathbb{R}$, but if a term $x$ had type $\mathbb{R}_{>0}$ then $x$ itself would \emph{not} strictly speaking have type $\mathbb{R}$; I stress again that every term has a \emph{unique} type. To make a term of type $\mathbb{R}_{>0}$ one has to give \emph{two} pieces of data: a real number, and a proof that it is positive. A term of type $\mathbb{R}_{>0}$ is an object corresponding to this pair, so strictly speaking it is not a real number, and a type theory based system will hold you to this. However of course there is a canonical map from $\mathbb{R}_{>0}$ to $\mathbb{R}$ -- you just throw away the proof. More generally, a type theory system could well have a \emph{coercion system}, consisting of a collection of ``invisible functions'' mapping types to other types in the way which mathematicians would expect. For example given a term of type $\mathbb{R}_{>0}$ it might well be possible to feed it into a function which is expecting a term of type $\mathbb{R}$; the system will just throw away the proof of positivity and use the underlying real number anyway. Mathematicians use these invisible functions everywhere, often without noticing. We have already mentioned above that in a foundational system the real numbers need to be built using one of the standard constructions, for example via Cauchy sequences. In particular a rational number is not literally a real number. However, taking Lean as an example, if one has a term \lstinline{x : ℚ} then one can simply write \lstinline{x : ℝ} to get the corresponding real number, although a careful inspection of the corresponding term will unearth the fact that the real number is actually called \lstinline{↑x}, indicating that a coercion has been applied. The coercion is a ring homomorphism, and Lean has a ``normalise casts'' tactic~\cite{normcast} which knows this and will apply theorems such as \lstinline{↑(x+y)=↑x+↑y} and \lstinline{↑(x*y)=↑x*↑y} automatically (before this tactic had been written, doing mathematics which involved switching between the naturals, integers and rationals could be quite frustrating because of these invisible maps). In summary then, type theory forces you to think more precisely about the actual objects you are working with, however tactics can be used to manipulate these objects the way we usually manipulate them. Learning how to ``steer'' mathematics in a theorem prover this way simply comes from practice.

\subsection{Inductive types}

I have already mentioned that in a type theory system the definition of the real numbers is the same as in a set theory system -- it is Cauchy sequences, or Dedekind cuts, or whatever your favourite construction of the reals is. Similarly the usual definitions of the rationals and integers as quotients work just as well in type theory as they do in set theory. But one place where the type-theoretic and set-theoretic foundations of mathematics differ is in the definition of the natural numbers. The natural numbers are a foundational object in mathematics -- they are typically the first example of an infinite object to be born -- so it is perhaps unsurprising that different foundational systems will treat them in different ways.

In ZFC set theory, the existence of the set of natural numbers is postulated as an axiom, namely the axiom of infinity. Type theories such as Lean's instead allow the user to define custom \emph{inductive types}. Such types include the naturals and other recursively-defined constructions. Implementation details of this so-called calculus of inductive constructions~\cite{CIC} differ between systems; the rest of this section explains details which are specific to Lean's type theory, but much of what I say applies to Coq and Agda, other popular type theory provers.

In Lean, the definition of the naturals looks like this:

\begin{lstlisting}
inductive nat
| zero : nat
| succ (n : nat) : nat
\end{lstlisting}

This definition says ``zero is a natural number, the successor of a natural number is a natural number, and that's it''. As one might guess, this inductive construction can be used to construct far more exotic types, but one can show that any type which can be defined using the rules of the calculus of inductive constructions corresponds to a set which can be built using the usual axioms of set theory.

Let us see what goes on under the hood when the naturals are defined as an inductive type. When such a definition is made, a new type \lstinline{nat} appears in the system, as does the term \lstinline{nat.zero} and the function \lstinline{nat.succ : nat → nat}. The latter terms are called \emph{constructors}: they are ways to make natural numbers. However one more thing also appears, namely the \emph{eliminator} for the type -- the object which enables the user to construct functions whose domain is the naturals and whose codomain is something else. It represents the idea that the only way that one can construct naturals is via \lstinline{nat.zero} and \lstinline{nat.succ}, and it states that to define a function out of the naturals, it suffices to (1) say where \lstinline{nat.zero} goes, and (2) to say where \lstinline{nat.succ n} goes, given where \lstinline{n} went. In other words, it is the principle of recursion.

So this is how new inductive types are born in Lean; after their definition they, together with their constructors and eliminator, are automatically added by the proof assistant to the system as new constants, or axioms, or however you would like to think of them. There are of course precise rules telling us the exact form of the eliminator for a given inductive type; we do not go into these here. From a foundational point of view this approach, where new axioms appear ``by magic''  as types are constructed, is very different to the set-theoretic viewpoint, however in~\cite{werner} it is shown that type theory with these constructions is equiconsistent with set theory. The strategy of the proof is to make a model of set theory within type theory, and to make a model of type theory within set theory. For a more precise statement, one has to be more precise about exactly what kind of type theory one is working with. For example Mario Carneiro's MSc thesis~\cite{mariothesis} shows that Lean's type theory is equiconsistent with ZFC plus countably many inaccessible cardinals.

It is worth noting, and quite amusing, that equality itself is defined as an inductive type in many type theory systems. This is in contrast to set theory, where equality is typically considered as part of the logic. Indeed, equality in type theory is generally more subtle than in set theory. Here is Lean's definition of equality:

\begin{lstlisting}
inductive eq {X : Type} : X → X → Prop
| refl (a : X) : eq a a
\end{lstlisting}

The slightly unnerving \lstinline{X → X → Prop}, bracketed as \lstinline{X → (X → Prop)}, means that equality is a function which takes in an element of \lstinline{X} and outputs a function which takes in an element of \lstinline{X} and outputs a Proposition, that is, a true-false statement. In other words, if \lstinline{a} and \lstinline{b} are terms of type \lstinline{X} then \lstinline{eq a b} is a true-false statement. Using the usual notation \lstinline{a = b} for \lstinline{eq a b}, we see that equality of terms of a type \lstinline{X} is an inductive type with one constructor, namely \lstinline{eq.refl a}, a proof that \lstinline{a = a}. It turns out that from this definition we can \emph{prove} all the usual properties of equality! The eliminator for the equality type is the \emph{substitution property}, that if $a=b$ then given a term of type $P(a)$ we can get a term of type $P(b)$. It is a rather pleasant game to go on from this to deduce that equality is both symmetric and transitive (for more details on this see for example~\cite{xena-equality}). Of course, whilst it is of interest to some to see how basic properties of equality can be proved within a type theory system, it is also worth stressing that to use a computer theorem prover one does not have to know anything about them.

\subsection{Dependent types}

Lean and Coq both use a version of type theory called dependent type theory, so it is perhaps worth taking some time to explain what a dependent type is. Imagine $X$ is a geometric object, for example a real manifold. Say that we have a vector bundle on $X$, that is, for each point $x$ of $X$ a vector space $V_x$ (which varies smoothly with $x$ in some appropriate sense). A section of this bundle is a function which takes as input a point $x$ in $X$ and outputs an element of $V_x$. From a foundational point of view there are two ways to think about such a section. One could regard this section as a function from $X$ to the disjoint union of the $V_x$, sending $x\in X$ to an element of $V_x$. Alternatively one could regard it as a slightly stranger kind of ``function'' which has domain $X$ but whose codomain varies according to the input. There are times in mathematics when taking the disjoint union of the codomains is a natural thing to do -- for example in the example above, the disjoint union of the $V_x$ is naturally a space $V$ sitting above~$X$. However there are also times when taking the disjoint union is quite unnatural. For example, in algebraic geometry one way of defining the sections of the structure sheaf on an affine scheme $\mathrm{Spec}(R)$ is functions which send a prime ideal~$P$ of~$R$ to an element of the localisation $R_P$ of $R$ at $P$, and the disjoint union of the $R_P$ as $P$ varies over the prime ideals of~$R$ has no natural algebraic structure. The set or type consisting of the disjoint union of these local rings is typically not part of the mental model which an algebraic geometer has when describing these sections.

These kinds of ``functions'' which have a well-defined domain, but a codomain which can vary according to the input, are called dependent functions. Not all proof assistants have such functions; for example Isabelle/HOL (a powerful proof assistant which contains a lot of analysis and analytic number theory) and various other HOL systems do not have them, which means that certain constructions in geometry are more convoluted than in Coq or Lean. See for example~\cite{isabelle-schemes}, which defines schemes in Isabelle/HOL but which has to build a new implementation of ring theory from scratch in order to do so.

\subsection{Examples}

Let us take a look at some examples of what mathematics looks like in a theorem prover based on type theory. I give these examples mainly to convince the reader who has been brought up using the language of set theory that there really is very little difference.

Here is what the claim that $\sqrt{2}$ is not rational looks like in Isabelle/HOL:

\begin{verbatim}
theorem sqrt2_not_rational:
  "sqrt 2 ∉ ℚ"
\end{verbatim}

You can see the proof on Isabelle's Wikipedia article~\cite{isabelle-wiki}. The fact that 2 is a term of a type and not a set, or an element of a set, is invisible.

Here is some more advanced mathematics, written in Coq:

\begin{verbatim}
Lemma prod_Cyclotomic n :
(n > 0)%N -> \prod_(d <- divisors n) 'Phi_d = 'X^n - 1.
\end{verbatim}

This is the statement that the product of the $d$th cyclotomic polynomials over $d\mid n$ is $X^n-1$. Note the hypothesis $n>0$, an assumption which a human would typically omit; computers are very picky with such ``edge cases''.

Here is the definition of a perfectoid ring in Lean, taken from the Lean perfectoid spaces website~\cite{perfectoid-site} which accompanies the article~\cite{bcmperfectoid}.

% \vfill\eject

\begin{lstlisting}
/-- A perfectoid ring is a Huber ring that is complete, uniform,
that has a pseudo-uniformizer whose p-th power divides p in the power bounded subring,
and such that Frobenius is a surjection on the reduction modulo p.-/
structure perfectoid_ring (R : Type) [Huber_ring R] extends Tate_ring R : Prop :=
(complete  : is_complete_hausdorff R)
(uniform   : is_uniform R)
(ramified  : ∃ ϖ : pseudo_uniformizer R, ϖ^p ∣ p in Rᵒ)
(Frobenius : surjective (Frob Rᵒ∕p))
\end{lstlisting}

The comment at the top of the code is the ``docstring'' for the code -- this is the human-readable explanation of what the Lean definition \lstinline{perfectoid_ring} represents, and this docstring is visible when you hover your cursor on the word \lstinline{perfectoid_ring} in some Lean code; if you are running the code in an IDE such as Microsoft VS Code then right-clicking on this word will jump you to the definition.

The Lean definition pretty much coincides with the human definition. If $R$ is a Huber ring which is a Tate ring (these are technical properties of topological rings), then we say $R$ is a perfectoid ring if it is complete, uniform and satisfies a couple of technical properties. The point to observe is that the computer code is no more or less difficult than the human definition.

\subsection{Foundations}

In my experience, mathematicians often have very little interest in the technicalities of the logical foundations of their subject -- they cannot list the axioms of set theory, but they know from experience what is ``legal mathematics''. The controversies of the early 20th century about whether nonconstructive methods are allowed in mathematical proofs have long ago died down; working mathematicians use the law of the excluded middle all over the place, and many use the axiom of choice in some form or another (indeed countable dependent choice can be invoked almost without one noticing). A typical research mathematician will have gone to at most one course on the foundations of mathematics; in such a course one typically learns that Zermelo--Frankel set theory with the axiom of choice, or ZFC, can be used as a foundation for much of mathematics. Indeed it can be used for essentially all of mathematics up until the 1960s; however Grothendieck's super-general cohomology theories developed in SGA4 introduced a new ``axiom of universes'' (the assertion that every set is an element of a set which is a model of ZFC). This axiom cannot be proved from the axioms of ZFC, by G\"odel's theorem. The original proofs of the Weil conjectures in theory used this axiom in the weak sense that at the time the only reference for \'etale cohomology was SGA4. However Deligne and others point out in SGA4$\frac{1}{2}$ that the theory of \'etale cohomology, and hence the proof of the Weil conjectures, can be set up within ZFC alone. Readers interested in the contortions that one has to go through in order to do this can look at the Set Theory section of the Stacks project, for example here~\cite[\href{https://stacks.math.columbia.edu/tag/000H}{Tag 000H}]{stacks-project}. For a more extreme example, see section~4 of~\cite{scholze-diamonds}, where we see a Fields Medallist forcing a more elaborate theory into ZFC.

My personal opinion is that whilst ZFC was a wonderful foundation for much of early 20th century mathematics, the lack of a universe axiom now means that it is becoming more and more of an effort to get parts of modern mathematics to fit into it. In books and papers dealing with infinity categories or condensed mathematics it is not at all uncommon to see universes showing up, and I do wonder whether now it is time for mathematicians to begin embracing universes, as Grothendieck was encouraging us to do since the 1960s. Coq's type theory and Lean's type theory both contain universes as part of the foundations, however mathematicians can choose not to use them if they so desire.

\section{The Future}

In this section I describe some of the plausible consequences of formalising mathematics in a computer theorem prover. I also highlight some things which I believe will remain out of reach for some time yet. Patrick Massot's more extensive observations~\cite{massot-note} are also well worth a read (indeed several of my ideas here were formed after conversations with Massot).

\subsection{A new kind of mathematical document}

Right now, an author of a textbook or research paper has to decide how much background material to assume, and which techniques they will regard as standard in the arguments they present. In other words, they have to decide where to start, and how fast to go. If a potential reader (for example a new PhD student, or an undergraduate interested in the area) does not have the necessary prerequisites then it will be far more difficult for them to get anything out of the paper.

Computer formalisation offers the possibility of a new kind of mathematical document, where the \emph{reader} can make the decisions about how much detail is visible. Patrick Massot has been experimenting with such documents. A preliminary version of his vision can be seen at his Sphere Eversion Project web pages~\cite{sphere-eversion}. This is a project whose main goal is to formalise in Lean a proof of Smale's theorem saying that a sphere can be turned inside out (or more formally, that there is a homotopy of immersions between the identity immersion of $S^2$ in $\mathbb{R}^3$ and the antipodal immersion). At the time of writing, the proof is not yet fully formalised, but it is only a matter of time. The blueprint is written in \LaTeX, but using plasTeX it has been converted into a web page with live Lean links. Right now these links take you to static web pages containing Lean code, but tools are currently being developed which will change this. Alectryon is a program available for Coq and Lean which can turn compiled code into web pages. Tools like Alectryon will enable us to make documents which will allow links to dynamic web pages displaying anything from mathematical details to interactive pictures, in a human-readable form, and which will allow one to keep digging right down to the axioms, although of course it is unlikely that anyone would like to go down this far.

There are already variants of this idea in existence, Lamport's idea of a ``structured proof'' came from a desire to encourage mathematicians to write far more details down in their papers, but one can see why such a proposal would not go down well. Here we can let automation do part of the work for us. The Metamath proof assistant also offers similar functionality already, because Metamath has very little automation and hence drilling down to the axioms is essentially the same as inspecting the proof.

One could also imagine error-free undergraduate textbooks also written in this way, where statements which a student cannot understand (perhaps because they are ambiguous) can be inspected in more details until difficulties are resolved.

\subsection{Semantic search in a mathematical database}

One thing that is not going to happen any time soon, is some kind of revolution where all mathematicians start writing all their papers in a formal proof assistant. Whilst one might expect a future where \emph{some} papers are partially, or even completely formalised in a theorem prover (see for example~\cite{gouezel-shchur}, \cite{strickland2019iterated} and~\cite{kjoshanssen2021parametrized}), this kind of approach will not become the norm any time soon. Faced with this reality, how will formalised mathematics be able to keep up with the frontiers of mathematics?

I have already mentioned Tom Hales' 2017 ``Big Conjectures'' talk at the Newton Institute in Cambridge. In the talk~\cite{hales-newt}, Hales argues for a formalised version of Math Reviews/Zentralblatt. That is, a website whose role is to formally \emph{state} the results being announced in the main mathematical journals. Note that such a project is nowhere near as far-fetched as the idea of formalising mathematical \emph{proofs} in real time; theorem \emph{statements} are far easier to formalise.

The issue with Hales' plan, as he points out in the talk, is that to be able to formalise statements of theorems in even a part of modern mathematics such as the Langlands philosophy, one would have to define all of the basic objects which mathematicians in this area use. In the Langlands philosophy this would include, but be by no means limited to, definitions of automorphic forms and automorphic representations, Galois representations, abelian varieties, the rings defined by Fontaine and used to do $p$-adic Hodge theory, schemes, all the cohomology theories used in the area, perfectoid spaces, adeles and ideles,\ldots. The Lean community has over the last few years pushed hard to get some of the main definitions of modern research mathematics into {\tt mathlib}. At Imperial College alone we currently have Oliver Nash developing the basics of the theory of Lie Algebras so we can talk about centres of universal enveloping algebras, Mar\'ia In\'es de Frutos-Fern\'andez developing the theory of adeles and ideles of global fields with an eye on the statements of class field theory, Amelia Livingston developing group and Galois cohomology, Jujian Zhang developing sheaf cohomology with an eye on GAGA, and Ashvni Narayanan developing the basics of Iwasawa theory in her PhD thesis. I have already mentioned the work of myself, Massot and Commelin defining perfectoid spaces. The work of Scott Morrison, Bhavik Mehta, Justus Springer and Adam Topaz has recently enabled us to start developing the theory of sheaves on sites and homological algebra, so cohomology theories are now not too far away. Of course much remains to be done, but we are hoping that the idea of being able to formally \emph{state} the theorems of Annals and Inventiones algebraic number theory papers in Lean will soon become a reality.

A related project is formalising tags in the Stacks Project. The Stacks Project~\cite{stacks-project} is a gigantic online database of algebraic geometry, freely accessible online. When printed out, it fills over 7000 pdf pages. Formalising all the proofs in the database would be an extremely arduous task involving many person-decades of work with current technology. In theory it is possible, however one would need a team who were experts in both algebraic geometry and in formalisation. Furthermore, for it to actually happen, the incentive structure in academic mathematics would have to change drastically. Publishing papers in prestigious computer science conference proceedings explaining how you developed the basic theory of Cohen--Macauley rings and modules in a theorem prover (and of course such work would be publishable in a prestigious computer science conference proceedings -- nobody has ever done it before) is perhaps not something which is recognised by promotions committees. 

However there is a solution available to us right now. Formalising just the \emph{definitions} and theorem \emph{statements} in the Stacks Project is a \emph{much} simpler task. Anybody interested in algebraic geometry would be more than welcome to learn Lean by attempting to formalise statements in Stacks Project tags. Point your web browser to the Lean Zulip instance~\cite{zulip} and ask where to get started in the {\tt \#new members} stream.

The reason that building such databases is important is that they will enable the community to build tools the likes of which mathematicians have never seen before. Let's imagine that all the definitions and theorem statements in the Stacks Project have been formalised in Lean or some other theorem prover. A ``hammer'' is a tool which runs inside a theorem prover and which can attempt to construct mathematical arguments by piecing together results in a database. The original hammer was Isabelle/HOL's \emph{Sledgehammer}~\cite{sledgehammer}. The cleverness behind such tools is the ability to isolate which of the many results in the database look the most useful, and to concentrate on these when attempting to prove the required result. Now consider a PhD student who is beginning to learn algebraic geometry. Such a student would then be able to ask the theorem prover a question, and the prover could attempt to use the database to answer the question positively (by piecing together a proof) or negatively (by producing a counterexample, like the website $\pi$-base~\cite{pi-base} is doing for counterexamples in topology). The resulting output of the computer would be able to explicitly point to references in the literature, or direct proofs of the claims it is making in its argument. This sort of tool -- computer assisted learning -- has the potential to beat the techniques currently used by PhD students (``google hopefully'', ``page through a textbook/paper hopefully'', ``ask on a maths website and then wait'', ``ask another human'') hands down. But as I have stressed before, the main thing which is missing is the database of theorems, and it is up to us to construct it. The sooner it is there, the sooner the tools will appear. And the bigger the database gets, the more powerful the tools will become.

\subsection{Checking proofs}

Some computer scientists have argued that mathematicians are sloppy, and our literature has errors in, and that this problem can be solved with computer proof assistants. Such an argument might initially look plausible, and I myself was a proponent of it a few years ago, but it does not stand up to much scrutiny. Firstly, the experts in our community know which results can be relied upon. Secondly, many errors are not serious and can be fixed. Thirdly, the more serious instances of this problem cannot be solved with computer proof assistants right now anyway. A great example is Mochizuki's claimed proof of the ABC conjecture~\cite{abc}. This proof has now been published in a serious research journal, however it is clear that it is not accepted by the mathematical community in general. One could challenge Mochizuki, or indeed anyone, to formalise the proof in a computer theorem prover. However this would be a completely unreasonable thing to do. A computer formalisation is not expected of other proofs appearing in our literature. Furthermore, the key sticking point right now is that the unbelievers argue that more details are needed in the proof of Corollary~3.12 in the main paper, and the state of the art right now is simply that one cannot begin to formalise this corollary without access to these details in some form (for example a paper proof containing far more information about the argument).

What \emph{would} however be feasible is for mathematicians to formalise \emph{parts} of technical work, or to get others to do so. There might be several reasons to do such a thing -- Commelin and his team have already shown that theorem provers can be used to check parts of complex proofs which humans might find it difficult to plough through, whilst learning about the mathematics in the process.

\subsection{Teaching}

I have heard students say ``I think my proof is OK'' when talking about their homework. Computer proof assistants are able to tell them immediately if this is so -- as long as the student has taken the trouble to learn the language of the proof assistant. Should we be teaching undergraduate mathematicians how to use computer proof assistants? I certainly think so. Patrick Massot in Orsay and myself at Imperial College London are both teaching undergraduate-level courses which do precisely this.

Students want feedback on their work as soon as possible. A computer proof assistant can supply it immediately.

Beginner students can be confused about the basics. What is the difference between $\forall \epsilon>0, \exists \delta>0, \ldots$ and $\exists \delta>0, \forall \epsilon>0, \ldots$? Once these systems become easier for mathematicians to use, students can experiment for themselves with well-chosen examples supplied by a lecturer and begin to understand what is going on. I was once told by a student ``I did not understand equivalence relations, so I formalised them in Lean, and then I understood them''. Forcing students to think pedantically and logically can be good for them.

It is however worth stressing that asking a weak student to both keep up with your course and to simultaneously learn how to use a computer theorem prover is clearly asking too much from that student. The provers need to become easier to use, perhaps with graphical interfaces and documentation more appropriate for mathematicians. Asking people to change the way they teach is of course asking a lot. However, mathematics education experts will be only too happy to tell us that our preferred medium -- ``write for an hour on a board'' -- is becoming less and less appropriate for our students, who like to learn things by watching 5 minute videos or playing with interactive toys. Can we make abstract mathematics more interactive? I suspect that we can. The more people who understand how to use these machines, the sooner the new ideas will come.

\subsection{Other ideas}

I do not claim to have exhausted the possibilities here. The people who designed the CD in the 1980s surely could not envisage music services like YouTube and Spotify, or the audiobook. The people who started to think about how to make typesetting of books look good on a computer screen surely did not envisage devices like the Kindle. It is time to look beyond how we usually teach and learn mathematics, and try to understand how we as a community of mathematicians can use the inevitable digitisation of mathematical material as a tool to make our lives, and the lives of our students, better. As Carneiro once said, you can't stop progress.

\bibliographystyle{amsalpha}

\bibliography{icm.bib}

\end{document}